\documentclass[11pt]{article}
\usepackage{amssymb}
\usepackage{graphicx}
\usepackage{latexsym}
\usepackage{amsmath,amscd}
\usepackage{bm}
\def\part#1{\frac{\partial\phantom{q}}{\partial#1}}

\newenvironment{rmk}{\begin{trivlist}\item[]{\bf Remark:} }
{\end{trivlist}}

\newenvironment{ex}{\begin{trivlist}\item[]{\bf Example:} }
{\end{trivlist}}
\newenvironment{exs}{\begin{trivlist}\item[]{\bf Examples:} }
{\end{trivlist}}
\newenvironment{prf}{\begin{trivlist}\item[]{\bf Proof:} }
{\hfill $\Box$ \end{trivlist}}

\newtheorem{thm}{Theorem}

\newcommand{\lie}[1]{\mathfrak{#1}}
\def\End{\mathop{\rm End}\nolimits}

\def\Hom{\mathop{\rm Hom}\nolimits}

\def\Ad{\mathop{\rm Ad}\nolimits}

\newcommand{\R}{\mathbf{R}}
\newcommand{\C}{\mathbf{C}}
\newcommand{\K}{\mathbf{H}}
\newcommand{\Z}{\mathbf{Z}}

\newcommand{\CP}{{\mathbf P}}

\begin{document}
\title{The hyperholomorphic line bundle }
 \author{Nigel Hitchin\\[5pt]}
 \maketitle
 \centerline{\it Dedicated to Klaus Hulek on the occasion of his 60th birthday}
 \vskip 1cm

\section{Introduction}
In a recent paper \cite{Hay}, A.Haydys introduced a natural line bundle with connection on a hyperk\"ahler manifold with an $S^1$-action of a certain type. The curvature is of type $(1,1)$ with respect to all complex structures in the hyperk\"ahler family and for this reason is called {\it hyperholomorphic}. In \cite{NJH} a description of this line bundle via  a holomorphic bundle on the twistor space was given and in this format calculated for a number of examples of interest to physicists. These are mostly moduli spaces of solutions to gauge-theoretic equations. 

In this article we give examples with a more geometrical flavour, in particular on minimal resolutions of Kleinian singularities and cotangent bundles of coadjoint orbits of a compact Lie group. We first approach the subject from the differential-geometric point of view, giving some explicit formulae, and then from the twistor viewpoint, where, as in \cite{NJH}, the holomorphic point of view  demonstrates  a  naturality which is not apparent from the explicit expressions. 

In a more general result, which contributes to the  examples, we show how the hyperholomorphic  bundle descends naturally in a hyperk\"ahler quotient, and  for the quotient by a linear action on flat space can be identified with a canonical hyperholomorphic  line bundle. 

\section{The differential geometric viewpoint}
\subsection{The hyperholomorphic connection}\label{defs}
Let  $M$ be  a hyperk\"ahler manifold with K\"ahler forms $\omega_1,\omega_2,\omega_3$ relative to complex structures $I,J,K$. If the de Rham cohomology class $[\omega_1/2\pi]\in H^2(M,\R)$  is in the image of the integral cohomology then there exists a  line bundle $L$ and  hermitian connection $\nabla$ with curvature $\omega_1$, unique up to tensoring with a flat $U(1)$ bundle. Since $\omega_1$ is of type $(1,1)$ with respect to the complex structure $I$, $L$ also has a holomorphic structure  defined by the $\bar\partial$-operator $\nabla^{0,1}$. Given a local holomorphic section $s$ of $L$, then $\omega_1=dd^c\log \Vert s\Vert^2/2$. Hence, if we multiply the hermitian metric by $e^{2f}$ the curvature of the connection on $L$ compatible with this new structure is 
$$F=\omega_1+dd^c\!f.$$
Suppose now we have a circle action which fixes $\omega_1$ but acts on the other forms by the transformation  $(\omega_2+i\omega_3)\mapsto e^{i\theta}(\omega_2+i\omega_3)$. The manifold $M$ must necessarily be noncompact for this. Suppose further that we have chosen a lift of the action to $L$. This implies  in particular the existence of a moment map  -- a function $\mu$ such that $i_X\omega_1=d\mu$ where $X$ is the vector field generated by the action. Then the result of Haydys \cite{Hay} (see also \cite{NJH}) is:

\begin{thm} The 2-form $\omega_1+dd^c\mu$ is of type $(1,1)$ with respect to complex structures $I,J,K$.
\end{thm}

Thus rescaling the natural metric by $e^{2\mu}$  gives a new connection which defines a holomorphic structure on $L$ relative to all complex structures in the quaternionic family. This is a hyperholomorphic connection, and $L$ is the hyperholomorphic bundle of the title. 

There are relatively few hyperk\"ahler metrics which one can write down explicitly but it is instructive to find the line bundle in these cases. 

\begin{ex}
Flat quaternionic space $\K^n$. Writing $\K^n=\C^n\oplus j\C^n$ we have 
$$\omega_1=\frac{i}{2}\sum_i(dz_i\wedge d\bar z_i+dw_i\wedge d\bar w_i),\quad \omega_2+i\omega_3=\sum_idz_i\wedge dw_i$$
and the action $(z,w)\mapsto (z,e^{i\theta} w)$ is of the required form. Then 
$$F=\omega_1+dd^c\mu=\frac{i}{2}\sum_i(dz_i\wedge d\bar z_i-dw_i\wedge d\bar w_i).$$
In the complex structure $I$ this is the trivial holomorphic line bundle with hermitian metric $h=(\Vert z\Vert^2-\Vert w\Vert^2)/2$.
\end{ex}

In the above we have specified a particular action of the circle  on the three K\"ahler forms $\omega_1,\omega_2,\omega_3$. More generally, if an irreducible  hyperk\"ahler manifold $M$ has a circle symmetry group then it acts on the three-dimensional space of covariant constant $2$-forms preserving the inner product. The action is either trivial, in which case it is called {\it triholomorphic}, or it leaves fixed a one-dimensional subspace  with an orthogonal complement on which the action is rotation by $n\theta$.  The case above is  $n=1$. This occurs for example on the cotangent bundle of a complex manifold where the action is scalar multiplication in a fibre and the symplectic form is the canonical one. 
 In the general case, $\Z_n\subset S^1$ preserves the three K\"ahler forms and so the quotient $M/\Z_n$ is a hyperk\"ahler orbifold with a circle action as above.  The local geometry of the hyperholomorphic bundle is then the same, but  the curvature form on $M$ is  $F=\omega_1+ndd^c\mu$.

In what follows we shall also consider  flat space as above but with the  action  $(z,w)\mapsto e^{i\theta}(z, w)$. Then $n=2$ since $(\omega_2+i\omega_3)\mapsto e^{2i\theta}(\omega_2+i\omega_3)$. The moment map $\mu=-(\Vert z \Vert^2+\Vert w\Vert^2)/2$ and so $F=\omega_1+2dd^c\mu=0$ and the hyperholomorphic line bundle is trivial as a line bundle with connection. This may seem uninteresting, but in Theorem \ref{quoti} we shall see how it defines the bundle for a hyperk\"ahler quotient of $\K^n$.

\subsection{Hermitian symmetric spaces}\label{herm}
Biquard and Gauduchon gave in \cite{BG} an explicit formula for a hyperk\"ahler metric which, in the complex structure $I$, is defined on the total space of the cotangent bundle of a hermitian symmetric space $G/H$. A circle action is given by multiplication of a cotangent vector by a unit complex number and  the form $\omega_2+i\omega_3$ is the canonical symplectic form on the cotangent bundle. 

If $p:T^*(G/H)\rightarrow G/H$ is the projection and $\omega$ is the K\"ahler form of the symmetric space  $G/H$ then the hyperk\"ahler metric is defined by
$\omega_1=p^*\omega+dd^ch$ where, for a cotangent vector $v$, $h$ is the quartic function on the fibres defined by $h(v)=(f(IR(Iv,v))v,v)$. Here $R(u,v)$ is the curvature tensor of $G/H$ and $f$ is the analytic function 
$$f(u)=\frac{1}{u}\left(\sqrt{1+u}-1-\log \frac{1+\sqrt{1+u}}{2}\right).$$
This function is applied to $IR(Iv,v)$ which is a non-negative hermitian transformation. In fact since the curvature of a symmetric space is constant we can also view the quadratic map $R(Iv,v)$ from $(\lie{g}/\lie{h})^*$ to $\lie{h}\subset \lie{g}$ as a multiple of the moment map for the isotropy action of $H$. 
The strange function $f(u)$ has the property that 
\begin{equation}
(uf(u))'=\frac{1}{2u}(\sqrt{1+u}-1)
\label{fu}
\end{equation}
We first calculate the moment map $\mu$ for the circle action. Since the action is purely in the fibres of the cotangent bundle we have
$$i_X\omega_1=i_X(p^*\omega+dd^ch)=i_Xdd^ch.$$
Now the action preserves both $h$ and the complex structure so $(di_X+i_Xd)d^ch={\mathcal L}_Xd^ch=d^c({\mathcal L}_Xh)=0$, which means that 
$i_X\omega_1=-d(i_Xd^ch)$ and we can take $\mu=-i_Xd^ch=(IX)(h)$.  The vector field $X$ was generated by  $v\mapsto e^{i\theta}v$ so $IX$ is generated by $v\mapsto \lambda^{-1} v$ for $\lambda\in \R^+$. 
Hence 
$$\mu(v)=\frac{\partial}{\partial \lambda}h(\lambda^{-1}v)\vert_{\lambda=1}.$$
But $h(v)=(f(u)v,v)$ where $u=IR(Iv,v)$ is homogeneous of degree $2$ in $v$ and so $\mu(v)=-2(uf'(u)v,v)-2(f(u)v,v)=-2((uf(u))'v,v)$. Using (\ref{fu}) we see that 
$$F=\omega_1+dd^c\mu=p^*\omega+dd^ck$$
where $k(v)=(g(IR(Iv,v))v,v)$ for the function $$g(u)=-\frac{1}{u}\left(\log \frac{1+\sqrt{1+u}}{2}\right).$$
This is an explicit formula for the curvature of the hyperholomorphic line bundle (assuming $\omega$ is normalized so that $[\omega/2\pi]$ is an integral class). 

Note that on the zero-section $v=0$, $F$ restricts to $\omega$ and is $S^1$-invariant. From \cite{Feix1},\cite{Feix2} we can say that this is the unique hyperholomorphic extension to $T^*(G/H)$ of this line bundle with connection on $G/H$. Later we shall view this in a more natural setting. 
\subsection{Multi-instanton metrics}\label{mult}
The most concrete examples of hyperk\"ahler metrics are the gravitational multi-instantons of Gibbons and Hawking \cite{GibHaw}. These are four-dimensional and in this dimension a hyperholomorphic connection is the same thing as an anti-self-dual one. The general Ansatz for this family of metrics consists of taking a harmonic function $V$ on an open set in $\R^3$, with its flat metric. Writing locally $\ast dV=d\alpha$ the metric has the form
$$g=V(dx_1^2+dx_2^2+dx_3^2)+V^{-1}(d\theta+\alpha)^2.$$
Then $\omega_1=Vdx_2\wedge dx_3+dx_1\wedge (d\theta +\alpha)$ is a K\"ahler form and similarly for $\omega_2,\omega_3$. 

An example is flat space $\C^2$ with a circle action $(z_1,z_2)\mapsto (e^{i\theta}z_1,e^{-i\theta}z_2)$. (Note that this action is triholomorphic,  and so is not of the type we have been considering). The quotient space is $\R^3$ with Euclidean coordinates $x_1=(\vert z_1\vert^2-\vert z_2\vert^2)/2, x_2+ix_3=z_1z_2$ and then the metric has the above form if $V=1/2r$. The flat space $\C^2\backslash\{0\}$ is here expressed as a principal circle bundle over $\R^3\backslash \{0\}$ and $d\theta+\alpha$ is the connection form for the horizontal distribution defined by metric orthogonality. The curvature of the connection is $d\alpha=\ast dV$ and the function $V^{-1/2}$ is the length of the vector field $Y$ generated by the action. 

The general case has the same principal bundle structure but the flat example shows that a $1/r$ singularity for $V$ does not produce a singularity in the metric: it is simply a fixed point of the circle action on the four-manifold. With this in mind, setting
$$V=\sum_{i=1}^{k+1}\frac{1}{\vert {\mathbf x}-\mathbf {a_i}\vert}$$
for distinct points $\mathbf {a_i}\in \R^3$ defines a nonsingular, complete hyperk\"ahler manifold $M$.

If the points ${\mathbf a_i}$ lie on the  $x_1$-axis then rotation about that axis induces an isometric circle action generating a vector field $X$. This involves lifting the action on $\R^3$ to the $S^1$-bundle with connection form $\alpha$, commuting with the circle action. Such a lifting is defined by a vector field of the form 
$X=X_H+fY$,
  where $X_H$ is the horizontal lift of $$x_2\frac{\partial}{\partial  x_3}-x_3\frac{\partial}{\partial x_2}$$ 
and, since $\ast dV$ is the curvature of the connection,  $i_X\!\ast \!dV=df$. Since ${\mathcal L}_XV=0$ the local existence of such an $f$ is assured. This means
$$df=(x_2V_2+x_3V_3)dx_1-x_2V_1dx_2-x_3V_1dx_3.$$
It follows that, with ${\mathbf a_i}=(a_i,0,0)$,  
\begin{equation}
f=\sum_{i=1}^{k+1}\frac{x_1-a_i}{\vert {\mathbf x}-\mathbf {a_i}\vert}+c.
\label{eff}
\end{equation}

The K\"ahler form $\omega_1$ is given by 
$\omega_1=Vdx_2\wedge dx_3+dx_1\wedge(d\theta+\alpha)$.
This is the curvature of a connection if its periods lie in $2\pi \Z$. Now the segment $[a_i,a_{i+1}]$  defines a one-parameter family of $S^1$-orbits which become single points over the end-points and  therefore form a 2-sphere in $M$. The manifold retracts onto a neighbourhood of a chain $[a_1,a_{2}],[a_2,a_{3}],\dots, [a_k,a_{k+1}]$ of $k$ such spheres which therefore generate $H_2(M,\Z)$. Integrating $\omega_1$ over the $i$th sphere gives $2\pi (a_{i+1}-a_i)$ and so for integrality   we require $a_{i+1}-a_i$ to be an integer. With these conditions we have, from Haydys's theorem,  a hyperholomorphic line bundle which, since the two actions commute, is invariant under the triholomorphic circle action on $M$.

 Kronheimer \cite{Kron0} showed that $S^1$-invariant instantons on the multi-instanton  space became monopoles on $\R^3$ with Dirac singularities at the marked points ${\mathbf a_i}$. Since the hyperholomorphic bundle is invariant we can define it this way by a $U(1)$ monopole:  a harmonic function $\phi$ on $\R^3$ and a connection $A$ such that $F=dA=\ast d\phi$. The Ansatz is 
\begin{equation}
\hat A=A-\phi V^{-1}(d\theta+\alpha)
\label{Ahat}
\end{equation}
where $\hat A$ is a local connection form on $M$. Thus
$$\omega_1+dd^c\mu = dA-d(\phi V^{-1})\wedge(d\theta+\alpha)-\phi V^{-1}\ast dV$$
and taking the interior product with $Y$ we obtain
$$i_Y(\omega_1+dd^c\mu)=-dx_1+i_Ydd^c\mu = d(\phi V^{-1}).$$
 Since $Y$ is triholomorphic, it preserves $I$ and since it commutes with $X$ it preserves $\mu$ so as in the previous section  $d(\phi V^{-1})=-dx_1-d(i_Yd^c\mu)$ and up to an additive constant,   
  $$\phi V^{-1}=-x_1-i_Yd^c\mu=-(x_1+i_Y(Ii_X\omega_1))=-(x_1-g(X,Y))$$
Now $g(X,Y)=V^{-1}f$. therefore  
$$\phi=-x_1V+f=-\sum_{i=1}^{k+1}\frac{a_i}{\vert {\mathbf x}-\mathbf {a_i}\vert}+c.$$

Note however that $A\mapsto A+c\alpha, \phi\mapsto \phi+cV$  takes $\hat A$ to $A+c\alpha-(\phi+cV)V^{-1}(d\theta+\alpha)=\hat A-c d\theta$ and so preserves the anti-self-dual curvature form $d\hat A$ (this absorbs the constant ambiguity too). We can therefore also take
$$\phi=\sum_{i=1}^k\frac{a_{k+1}-a_i}{\vert {\mathbf x}-\mathbf {a_i}\vert}+c$$
and since the coefficients $a_{k+1}-a_i$ are integers, this is a genuine $U(1)$-monopole which satisfies the Dirac quantization condition. 

There remains the question of the constant $c$. This is not in general zero since $k=1$ is flat space and we have seen the non-zero curvature of the connection in the previous section. Here we have by construction also a circle action which preserves all three K\"ahler forms so given one lifting of the rotation action on $\R^3$ to $M$ we can compose with a homomorphism to the triholomorphic circle to obtain another. The constant $c$ will then change by  $2\pi n, n\in \Z$.  

\begin{rmk} When $c=0$ the curvature $F$ is a linear combination of   ${\mathcal L}^2$ harmonic forms \cite{Rub},\cite{HHM}. In this case if $k=2m$ and ${\bf x}$ lies on the $x_1$-axis with $a_m\le x_1\le a_{m+1}$ then (\ref{eff}) shows that $f=0$. Note for future reference that this  means that the  middle $2$-sphere in the chain is point-wise fixed by the circle action.
\end{rmk}

The complex structure $I$ for the metrics above  is the minimal resolution of the Kleinian singularity $xy=z^{k+1}$. There are, thanks to Kronheimer \cite{Kron1}, hyperk\"ahler metrics on all such resolutions. These are produced by a finite-dimensional hyperk\"ahler quotient construction and this is semi-explicit -- the quotient metric of a subspace of flat space defined by a finite number of quadratic equations -- but the hyperholomorphic line bundle is well adapted to the quotient construction.
\subsection{Hyperk\"ahler quotients}
The hyperk\"ahler quotient construction of \cite{HKLR} proceeds as follows. Given a hyperk\"ahler manifold with a triholomorphic action of a Lie group $G$ we have, under appropriate conditions, three moment maps $\nu_1,\nu_2,\nu_3$ corresponding to the three K\"ahler forms $\omega_1,\omega_2,\omega_3$ and hence a vector-valued moment map $\nu:M\rightarrow \lie{g}^*\otimes \R^3$. Then, assuming $G$ acts freely on $\nu^{-1}(0)$, the manifold $\bar M=\nu^{-1}(0)/G$ with its quotient metric is hyperk\"ahler. 

In our situation we have a distinguished complex structure $I$ preserved by a circle action. The construction can then be viewed in a slightly different way. Firstly $\nu_c=\nu_2+i\nu_3$ is holomorphic with respect to $I$ and so the zero set $M_0=\nu_c^{-1}(0)$ is a complex submanifold of $M$ and hence $\omega_1$ restricts to it as a K\"ahler form. The group $G$ preserves $M_0$ and $\nu_3$ is the moment map for the restriction of $\omega_1$. Hence the hyperk\"ahler quotient is the symplectic quotient of $M_0$ by this action. 

\begin{thm} \label{quot} Suppose $M$ has a circle action as in Section \ref{defs}, commuting with $G$, so that the hyperk\"ahler quotient $\bar M$ has an induced action. Then the hyperholomorphic line bundle on $M$ descends naturally to the hyperholomorphic line bundle of $\bar M$. 
\end{thm}

\begin{prf} First recall that for a symplectic manifold $(N,\omega)$ with $[\omega/2\pi]$ integral there is a line bundle -- the prequantum line bundle -- with a unitary connection whose curvature is $\omega$. Given a lift of the action of a group $G$, the invariant sections on the zero set of the moment map define the prequantum line bundle on the symplectic quotient. 

To see this more concretely, let $Y$ be the vertical vector field of the principal $U(1)$-bundle $P$,  $X_a$  the vector field on $N$ given by $a\in \lie{g}$ and $\mu:N\rightarrow \lie{g}^*$ the moment map. Then a lift commuting with the $U(1)$-action is defined by $(X_a)_H+\langle \mu,a\rangle Y$ where $(X_a)_H$ is the horizontal lift. 
 An arbitrary   section of the line bundle is defined by a function $f$ on $P$, equivariant under the vertical action,  and an invariant section  satisfies $((X_a)_H+\langle \mu,a\rangle Y)f=0$. Thus on $\mu^{-1}(0)$ we have  $(X_a)_Hf=0$  which means that the section is covariant constant along the $G$-orbits. Hence the connection is pulled back from the symplectic quotient $\mu^{-1}(0)/G$. 

This is the construction for a symplectic manifold. Now suppose we take our hyperk\"ahler manifold with circle action and commuting triholomorphic $G$-action with hyperk\"ahler moment map $\nu$. We want to apply the above to $N=M_0=\nu_c^{-1}(0)$ for the symplectic quotient of $M_0$ is the hyperk\"ahler quotient of $M$. Now the circle action does not preserve $\omega_2+i\omega_3$ but it acts on $d\nu_c=d(\nu_2+i\nu_3)$ by multiplication by $e^{i\theta}$.  If we make a choice of moment map so that the action on $\nu_c$ is the same scalar multiplication, then the action will preserve $M_0=\nu_c^{-1}(0)$. Moreover, $\mu$ restricted to $M_0$, is the moment map for $\omega_1$  restricted to $M_0$.

The  line bundle with hyperholomorphic connection on $M$, and hence its restriction to $M_0$, was obtained from the prequantum line bundle by rescaling the hermitian metric by $e^{2\mu}$.  By what we have just seen, this  descends to $\bar M$, the symplectic quotient of  $N=M_0$. However, $G$ commutes with the circle action and so $\mu$ is $G$-invariant. It is also the  moment map for the induced action on the quotient, and 
 it follows that  rescaling the prequantum hermitian metric on $\bar M$ gives the hyperholomorphic bundle. 
\end{prf}

One other aspect of the quotient is that it comes equipped with a canonical principal $G$-bundle with a hyperholomorphic connection. Indeed $\nu^{-1}(0)/G=\bar M$ and $\nu^{-1}(0)$ is the total space of the principal $G$-bundle. The induced metric defines an orthogonal subspace in the tangent space to the  orbit directions and this is the horizontal space of a connection, which is hyperholomorphic. A differential-geometric proof of this was give in  \cite{GN} but it can be seen very naturally from the twistor space point of view  which we carry out in the next section. In fact, with fewer formulae and more geometry, the hyperholomorphic bundle appears much more naturally using holomorphic techniques. 
\section{The twistor viewpoint}
\subsection{The holomorphic bundle}
This section is essentially a review of the construction in \cite{NJH}. 
The twistor space $Z$ of a hyperk\"ahler manifold $M$ is the product $Z=M\times S^2$ given the complex structure $(I_{\bf u},I)$ where $I_{{\bf u}}=u_1I+u_2J+u_2K$ for a unit vector ${\mathbf u}\in \R^3$ and where the second factor is the complex structure of $S^2=\CP^1$. The projection $\pi:Z\rightarrow \CP^1$ is holomorphic and for each $x\in M$, $(x,S^2)$ is a holomorphic section, a {\it twistor line}. 

The fibre over ${\bf u}\in S^2$ is the hyperk\"ahler manifold $M$ with complex structure defined by ${\bf u}$ but it also has a holomorphic symplectic form relative to this complex structure. Using an affine coordinate $\zeta$ on $\CP^1$ where $u_2+iu_3=2\zeta/(1+\vert \zeta\vert^2)$ the complex structures $ I, -I$ are defined by $\zeta=0,\infty$ and the holomorphic symplectic form is $(\omega_2+i\omega_3)+2i\omega_1\zeta+(\omega_2-i\omega_3)\zeta^2$. Globally, this is a twisted relative 2-form $\omega$: a holomorphic section of $\Lambda^2 T^*_{Z/\CP^1}(2)$ where the $(2)$ denotes the tensor product with the line bundle $\pi^*{\mathcal O}(2)$, reflecting the quadratic dependence on $\zeta$. 
\begin{ex}
The twistor space for flat $\K^n$ is the total space of the vector bundle $\C^{2n}(1)$ over $\CP^1$. This is given by holomorphic coordinates $(v,\xi,\zeta)\in \C^{2n+1}$ on the open set $U$ defined by $\zeta\ne \infty$ and $(\tilde v,\tilde \xi,\tilde \zeta)$ for $V$  by $\zeta\ne 0$, with identification $(\tilde v,\tilde \xi,\tilde \zeta)=(v/\zeta,\xi/\zeta,1/\zeta)$ over $\zeta\in \C^*$. In these coordinates $Z$ is expressed as a $C^{\infty}$-product by the map $(z,w,\zeta)\mapsto (z+\zeta\bar w,w-\zeta \bar z,\zeta)$. 
\end{ex}
If a bundle on $M$ has a hyperholomorphic connection its  curvature is of type $(1,1)$ with respect to all complex structures parametrized by $\zeta$ and it follows that  its pull-back to $Z=M\times S^2$ has a holomorphic structure. Conversely any holomorphic vector bundle on $Z$ which is trivial on the twistor lines $(x,S^2)$  defines a hyperholomorphic connection on a vector bundle over $M$. This is the hyperk\"ahler version of the Atiyah-Ward  result for anti-self-dual connections. For a line bundle the triviality on twistor lines is simply the vanishing of the first Chern class. To get a unitary connection we impose a reality condition. It follows that  to describe a hyperholomorphic line bundle on $M$ we simply look for a holomorphic line bundle $L_Z$ on $Z$ determined by the circle action.  

\begin{ex}
In flat space with the action $(z,w)\mapsto (z,e^{i\theta}w)$ one can calculate the line bundle directly. The $(1,0)$-forms on $Z$ for $\zeta\ne \infty$ are spanned by $dz_i+\zeta d\bar w_i,dw_i-\zeta d\bar z_i, d\zeta$ and then with 
$$\log h_U=\frac{1}{2}\sum_i z_i\bar z_i-w_i\bar w_i+\zeta\bar z_i\bar w_i+\bar\zeta z_i w_i$$ 
we find
$$\bar\partial_Z \log h_U= \frac{1}{2}\sum z_iw_id\bar\zeta+z_id\bar z_i-w_id\bar w_i+\bar\zeta d(z_iw_i)$$
and hence $\bar\partial_Z \partial_Z \log h_U=(\sum _i-dz_id\bar z_i+dw_id\bar w_i)/2$, the curvature of the hyperholomorphic line bundle, on the open set $U$. Defining $\log h_V=-\log h_U(-1/\bar\zeta)$ on $V$, the pair $(h_U,h_V)$  defines a hermitian metric on the line bundle with holomorphic transition function on $U\cap V$ 
$$\exp(-\sum_iv_i\xi_i/2\zeta).$$
\end{ex}
The link between the differential geometric and holomorphic points of view is proved in \cite{NJH}. In fact  the line bundle $L_Z$  is essentially the prequantum line bundle for the family of holomorphic symplectic manifolds defined by $Z$. 

To understand this, and to see where the circle action enters in the construction, first note that since $\omega_2+i\omega_3$ transforms as $(\omega_2+i\omega_3)\mapsto e^{i\theta}(\omega_2+i\omega_3)$, differentiating with respect to $\theta$ we have $\omega_2={\mathcal L}_X\omega_3= di_X\omega_3$ and so $\omega_2$ and similarly $\omega_3$ are exact. Thus the $2$-form $(\omega_2+i\omega_3)/2i\zeta+\omega_1+(\omega_2-i\omega_3)\zeta/2i$ has the same cohomology class as $\omega_1$ for any $\zeta$ and is therefore, given the integrality condition on $[\omega_1/2\pi]$, the curvature of a line bundle on $M$. In the complex structure at $\zeta$, $(\omega_2+i\omega_3)/2i\zeta+\omega_1+(\omega_2-i\omega_3)\zeta/2i$ is of type $(2,0)$ therefore has no $(0,2)$ part:  hence the bundle has a holomorphic structure.

Now observe that the induced circle action on the twistor space generates a holomorphic vector field $V$ on $Z$. Since the action fixes $\pm I$, $V$ projects to the vector field $i\zeta d/d\zeta$ on $\CP^1$ vanishing at $\zeta=0,\infty$. This is a holomorphic section $s$ of ${\mathcal O}(2)$ and so the $2$-form we wrote above, $(\omega_2+i\omega_3)/2i\zeta+\omega_1+(\omega_2-i\omega_3)\zeta/2i$ is, on a specific fibre,  the restriction of the meromorphic relative differential form $\omega/2is\in \Omega^2_{Z/{\CP^1}}$. It turns out that this relative form is the restriction of a closed meromorphic 2-form ${ F_Z}$ on $Z$, which is the curvature of a meromorphic connection on the holomorphic line bundle.

\begin{thm} \cite{NJH} \label{preq} The line bundle $L_Z$ on the twistor space  $Z$  admits a meromorphic connection such that
\begin{itemize}
\item
there are simple poles at $\zeta=0,\infty$
\item
the curvature ${ F_Z}$ restricts to 
$$\frac{1}{2i\zeta}(\omega_2+i\omega_3)+\omega_1+\frac{1}{2i}(\omega_2-i\omega_3)\zeta$$ on each fibre over $\C^*\subset \CP^1$
\item
$i_V{ F_Z}=0$ where $V$ is the vector field generated by the circle action.
\end{itemize}
\end{thm}
\begin{rmk} Suppose the holomorphic vector field $V$ integrates to a $\C^*$-action. Then as ${ F_Z}$ is closed, the last property  tells us that this action gives a symplectic isomorphism between any of the holomorphic symplectic manifolds over $\zeta\in \C^*$.
\end{rmk}
Given that such a connection exists, the line bundle is essentially  uniquely determined by the residue of the connection, for given any two such bundles $L,L'$ with the connections as above and with the same residue at $\zeta=0,\infty$, the resulting   holomorphic connection on $L'L^*$ would have a  curvature which is a holomorphic $2$-form. But the normal bundle of a twistor line is $\C^{2n}(1)$ and so $T^*_Z\cong \C^{2n}(-1)\oplus {\mathcal O}(-2)$ on such a line. It follows that there are no holomorphic forms of positive degree on a twistor space since there is a twistor line through each point. Hence the connection on  $L'L^*$ is flat and this is in any case the ambiguity in choosing a prequantum connection. 

The residue is  canonically determined by the data of the action as follows (see \cite{NJH} for details). Since the connection has a singularity on a divisor of ${\mathcal O}(2)$, its residue  will be a section of $T^*_Z(2)$ on that divisor. Now since $T_{\CP^1}\cong{\mathcal O}(2)$ the projection $\pi:Z\rightarrow \CP^1$ gives an exact sequence of bundles:
$$0\rightarrow {\mathcal O}\rightarrow T^*_Z(2)\rightarrow T^*_{Z/\CP^1}(2)\rightarrow 0$$
 and the twisted relative form $\omega$ identifies $T_{Z/\CP^1}$ with $T^*_{Z/\CP^1}(2)$. The resulting extension
 $$0\rightarrow {\mathcal O}\rightarrow E\rightarrow T_{Z/\CP^1}\rightarrow 0$$
 can be identified with $TP/\C^*$ where $P$ is the holomorphic principal bundle of the prequantum line bundle for the real symplectic form $\omega_1$. The vector field $V$ on $Z$ is tangential to the fibres at $\zeta=0,\infty$ and is $X$ itself. The moment map defines an invariant lift to $P$ and hence a section of $TP/\C^*$. Under the isomorphism above, this is the residue of the connection. If we restrict it as a form to the fibre $\zeta=0$ it is $i_V(\omega_2+i\omega_3)/2i$
 
\begin{exs} 

\noindent i) For flat space with the action $(z,w)\mapsto e^{i\theta}(z,w)$ the line bundle $L_Z$ is trivial and the connection with the trivial action is just the meromorphic one-form
$$\frac{1}{2\tilde\zeta}\sum_i \tilde \xi_id\tilde v_i-\tilde v_id\tilde \xi_i=\frac{\zeta}{2}\sum_i \frac{\xi_i}{\zeta}d\frac{ v_i}{\zeta}-\frac{ v_i}{\zeta}d\frac{ \xi_i}{\zeta}=\frac{1}{2\zeta}\sum_i  \xi_id v_i- v_id \xi_i.$$
With the action $u\mapsto e^{in\theta} u$ it is 
\begin{equation}
2\pi i n \frac{d\zeta}{\zeta}+\frac{1}{2\zeta}\sum_i  \xi_id v_i- v_id \xi_i
\label{conn}
\end{equation}
\noindent ii) Flat space with the other action  $(z,w)\mapsto(z, e^{i\theta}w)$ requires  local connection forms $A_U,A_V$ such that $A_V=A_U+g_{UV}^{-1}dg_{UV}$. Define  
 $$A_{U}=\frac{1}{2\zeta}\sum v_id\xi_i\qquad A_V=-\frac{1}{2\tilde \zeta}\sum_i\tilde v_id\tilde \xi_i$$
 then on $U\cap V$ 
 $$A_V-A_U=-\frac{\zeta}{2}\sum_i\frac{\xi_i}{\zeta}d\frac{v_i}{\zeta}-\frac{1}{2\zeta}\sum_iv_id\xi_i=-d\left(\frac{1}{2\zeta}\sum_iv_i\xi_i\right).$$
 \end{exs}
 \subsection{Hyperk\"ahler quotients} 
 In the twistor formalism the hyperk\"ahler quotient is a very natural operation: it is just the fibrewise holomorphic symplectic quotient as long as the holomorphic vector fields generated by $G$ integrate  to an action of the  complexification $G^c$.  Each $a\in \lie{g}$ gives a holomorphic vector field $V_a$ tangential to the fibres of $\pi:Z\rightarrow \CP^1$ and the three moment maps for $V_a, a\in \lie{g}$ give a complex section 
 $${\bm \nu}=(\nu_2+i\nu_3)+2i\nu_1\zeta+(\nu_2-i\nu_3)\zeta^2$$
 of $\lie{g}^c(2)$. The twistor space $\bar Z$ of the hyperk\"ahler quotient is then simply ${\bm \nu}^{-1}(0)/G^c$ where the metric plays a role in determining the stable points for this quotient by a complex group. With this viewpoint the descent of the hyperholomorphic bundle through a quotient is, given Theorem \ref{preq}, the descent of the prequantum line bundle in a symplectic quotient (it is straightforward to check that the residue descends appropriately). 
 
 As we saw in the previous section, a hyperk\"ahler  quotient brings with it a canonical  hyperholomorphic $G$-bundle.  In fact, in the twistor interpretation, ${\bm\nu}^{-1}(0)$ is a principal $G^c$-bundle over the twistor space $\bar Z={\bm\nu}^{-1}(0)/G^c$ and it satisfies the reality condition to  define a hyperholomorphic principal $G$-bundle over $\bar M$.  A homomorphism $G\rightarrow U(1)$ then defines a hyperholomorphic line bundle and this raises the obvious question about whether, given a circle action, this is the hyperholomorphic bundle of the title.  
 
 In fact for a manifold to be a smooth hyperk\"ahler quotient of flat space such homomorphisms must exist. The standard moment map for a linear action  is quadratic and the origin lies in ${\bm\nu}^{-1}(0)$, so for smoothness we must change this by a constant. Equivariance however demands that the constant   is an invariant in $\lie{g}^*$: a homomorphism from $\lie{g}$ to $\R$. 
 
 Consider flat space $\K^n$ as a right $\K$-module, then $U(n)\subset Sp(n)$ is the subgroup commuting with left multiplication by  $e^{i\theta}$: this is  a distinguished complex structure $I$. Let $G\subset U(n)$ and $c\in \lie{g}^*$ be a $G$-invariant element. If $c$ is integral it corresponds to a homomorphism $\chi:G\rightarrow U(1)$. Let $\nu$ be the  standard quadratic hyperk\"ahler moment map for the  linear action, then taking the reduction at $\nu=(c,0,0)$,  the cohomology class of the K\"ahler form $\omega_1$ lies in $2\pi H^2(\bar M,\Z)$. Indeed the integrality for $c$ gives a lift of the $G$-action to the prequantum line bundle on $\nu_c^{-1}(0)$ which descends. 
 
 \begin{thm}\label{quoti} If the hyperk\"ahler quotient $\bar M$ of $\K^n$ by $G$ with $\nu=(c,0,0)$ is smooth, then the hyperholomorphic line bundle is $\nu^{-1}(c,0,0)\times_G \C$ endowed with the canonical connection,
where $G$ acts via $\chi:G\rightarrow U(1)$.   \end{thm} 

\begin{prf} From the twistor point of view the line bundle $L_Z$ on the quotient is defined by the property that local sections are the same as local $G^c$-invariant sections of the holomorphic line bundle on ${\bm\nu}^{-1}(0)$.  For flat space and the circle action above the latter, as we observed in Section \ref{defs},  is a trivial holomorphic bundle but has a non-trivial action defined by $\chi$. Thus on $Z$ the line bundle is associated to the principal $G^c$-bundle  ${\bm\nu}^{-1}(0)$ by $\chi$.
\end{prf}

\begin{exs}

\noindent i) The simplest example is the cotangent bundle of a complex Grassmannian, one of the Hermitian symmetric spaces of Section \ref{herm}. In this case the flat space is $M=V\oplus jV$ for $V=\Hom(\C^k,\C^n)$ and $G=U(n)$ acting in the obvious way. There is just a one-dimensional space of of invariant elements in $\lie{g}^*$ and $H^2(\bar M,\Z)\cong \Z$. Notice that $-1$ acting on the vector space is represented by $-1\in U(n)$  and hence acts trivially on the quotient. It is therefore $e^{2i\theta}$ which acts effectively on the quotient. Since $e^{i\theta}$ acts on $\omega_2+i\omega_3$  in flat space by multiplication by $e^{2i\theta}$, on the quotient the induced action is the standard one: in fact the fibre action on the cotangent bundle. 

\noindent ii) Taking $M=V\oplus jV$ where $V=\End \C^k\oplus \Hom(\C^k,\C^n)$ and $G=U(k)$ one obtains the moduli space of $U(n)$-instantons on $\R^4$ of charge $k$ or, with a non-zero moment map, the moduli space of noncommutative instantons. For $n=1$ this is the Hilbert scheme $(\C^2)^{[k]}$ of $k$ points on $\C^2$ and  the hyperholomorphic line bundle with complex structure $I$ is defined by the exceptional divisor. The  circle action is induced from scalar multiplication on $\C^2$ and so  the action on $\omega_2+i\omega_3$ is multiplication by $e^{2i\theta}$, since on the open set of $(\C^2)^{[k]}$ consisting of the configuration space of $\C^2$  the symplectic form is the sum $k$ copies of  $dz\wedge dw$.

\noindent iii)  In \cite{Kron1} Kronheimer constructed asymptotically locally Euclidean hyperk\"ahler metrics on minimal resolutions of Kleinian singularities ($\C^2/\Gamma$ for $\Gamma\subset SU(2)$ a finite group) by the quotient construction.  The construction is as follows. 
  Let $R=L^2(\Gamma)$ be the regular representation, $\C^2$ the basic representation from $\Gamma\subset SU(2)$ and put $M=(\C^2\otimes \End(R))^{\Gamma}$. Since $\End(R)$ has real structure $A\mapsto A^*$ and $SU(2)\cong Sp(1)$ this is a quaternionic vector space and  the group  $G=U(R)^{\Gamma}$ acts  as quaternionic unitary transformations. The ALE space appears as a hyperk\"ahler quotient of $M$ by the action of $G$. If $R_0,\dots,R_k$ are the irreducible representations of $\Gamma$, of dimension $d_i$ then 
 $$R=\bigoplus_i \C^{d_i}\otimes R_i$$ and so $U(R)^{\Gamma}\cong U(d_0)\times \dots \times U(d_k)$. From the McKay correspondence each $R_i$ labels a vertex of an extended Dynkin diagram of type $A,D,E$ and then 
 \begin{equation}
 M=\bigoplus_{i\rightarrow j} \Hom(\C^{d_i},\C^{d_j})
 \label{dyn}
 \end{equation}
 the sum taken over the edges of the diagram, once with each orientation. As shown in \cite{Kron1}, the invariant subspace of $\lie{g}^*$ can be identified with the Cartan subalgebra of the Lie algebra of type $A,D,E$ as can the cohomology $H^2(\bar M,\R)$, with $H_2(\bar M,\Z)$ the root lattice.  The case of $A_k$ is the multi-instanton metric of Section \ref{mult}, where the chain of $2$-spheres constructed explicitly realizes the Dynkin diagram of type $A_k$.
 
Here the circle action on the symplectic form of the quotient  will be the standard one if there is an element in $G$ which acts as $-1$.  For this, from (\ref{dyn}) we need to show that there exist $c_i= \pm 1$, $0\le i\le k$,  such that if $i,j$ are joined by an edge of the extended Dynkin diagram then $c_ic_j=-1$. For $A_1$ this is trivial. Consider the extended Dynkin diagram (for $k>1$) as a simplicial complex, then  this is the same as asking that the $\Z_2$-cocycle associating $-1$ to each $1$-simplex is a coboundary. The diagrams of type $D_k,E_6,E_7,E_8$ are contractible and so have zero first cohomology so this is true. For $A_k$ the diagram is homeomorphic to a circle and the  cohomology class in $H^1$ vanishes if there is an even number of  edges, which is when $k$ is odd. 

Now the $D,E$ Dynkin diagrams have a trivalent vertex which, in the presence of our circle action, corresponds  to a rational curve of self-intersection $-2$ which is pointwise fixed, since there cannot  be just three fixed points. And, as pointed out in Section \ref{mult}, when $k$ is odd, the central curve in the $A_k$ case is fixed.

 In these cases, with respect to the complex structure $I$, we have a rational curve of self-intersection $-2$ and a neighbourhood of  such a  curve is biholomorphic  to a neighbourhood of the zero section of the cotangent bundle. Moreover the  circle action is the standard scalar multiplication in the fibre.  Applying \cite{Feix1} this means that the Kronheimer metric with circular symmetry is the unique hyperk\"ahler extension of the induced metric on the distinguished $2$-sphere.

 \end{exs}

\subsection{Coadjoint orbits} 
The Hermitian symmetric spaces which we considered in Section \ref{herm} are special cases of coadjoint orbits of compact semi-simple Lie groups with their canonical K\"ahler structure. There is a very natural description of the twistor space of a hyperk\"ahler metric on the cotangent bundle of such a space  due originally to Burns \cite{Burns}. In fact, that paper only asserts the existence of such a metric in a neighbourhood of the zero section, but it was written before hyperk\"ahler quotients, and in particular the infinite-dimensional gauge-theoretic versions, were discovered. Much later, armed with a knowledge of existence theorems for Nahm's equations, Biquard \cite{Biq} revisited this description being assured  of global existence. The action of scalar multiplication in the cotangent fibres by $e^{i\theta}$ gives a circle action and we shall now seek a concrete description  of the line bundle $L_Z$ using Burns's approach. 

Let $G$ be a semisimple compact Lie  group with a bi-invariant metric and $z\in \lie{g}$ be an element with centralizer $H$. Then in the complex group $G^c$ there are parabolic subgroups $P_+,P_-$ with $P_+\cap P_- = H^c$. The real (co)adjoint orbit $G/H\cong G^c/P_+\cong G^c/P_-$, and the complex coadjoint orbit is $G^c/H^c$.

The Lie algebra $\lie{p}_+=\lie{h}+ \lie{n}_+$ where $\lie{n}_+$ is nilpotent,   $z\in \lie{h}$ by definition and we define two complex manifolds
$$Z_0=G^c\times_{P_+}\{\C\cdot z+ \lie{n}_+\}\qquad
Z_{\infty}=G^c\times_{P_-}\{\C\cdot z+ \lie{n}_-\}.$$
Since $P_{\pm }$ fixes $z$ modulo ${\mathbf n}_{\pm}$, the coefficient of $z$ defines a projection $\pi_0:Z_0\rightarrow \C$ and similarly for $Z_{\infty}$. The fibre over $0$ is the cotangent bundle  $T^*(G^c/P_+)\cong G^c\times_{P_+} \lie{n}_+$ and for $\zeta\ne 0$, $G^c\times_{P_+}\{\zeta z+ \lie{n}_+\}$ is an affine bundle over $G^c/P_+$. 

There is another description, however,  for  $(g,\zeta z+x_+)\mapsto( \Ad g(\zeta z+x_+),\zeta)$ identifies the fibre at $\zeta\ne 0$ with the $G^c$-orbit of $\zeta z$. Note that the map $z\mapsto \zeta z$ defines an isomorphism with the orbit of $z$ which  is not symplectic for the canonical Kostant-Kirillov form $\omega_{\mathrm{can}}$ but is for its rescaling $\omega_{\mathrm{can}}/\zeta$. 

The twistor space is obtained by  identifying  $Z_0,Z_{\infty}$ over $\zeta\in \C^*$ by $(x,\zeta)\mapsto (\zeta^{-2}x,\zeta^{-1})$. Then the two projections define $\pi:Z\rightarrow \CP^1$ and $\omega_{\mathrm{can}}$ defines the twisted relative symplectic form. 

We define line bundles $L_+,L_-$ over $Z_0,Z_{\infty}$ by pulling back the prequantum line bundle on $G/H=G^c/P_{\pm}$ using the projections  $p_0:Z_0\rightarrow G^c/P_+,  p_{\infty}: Z_{\infty}\rightarrow G^c/P_-$. Then to define a line bundle $L_Z$ on $Z$ we need an isomorphism between $L_+$ and $L_-$ over $\C^*\subset \CP^1$. But the prequantum line bundle is homogeneous, defined by representations $\chi_{\pm}: P_{\pm}\rightarrow \C^*$, and these agree on $H^c=P_+\cap P_-$. This therefore gives an isomorphism  $p_+^*L_+\cong p_-^*L_-$ on $Z_0\cap Z_{\infty}\cong G^c/H^c\times \C^*$.

To show that this truly is the twistor version of the hyperholomorphic bundle we may simply note that it does generate a hyperholomorphic line bundle but by the invariance of the construction it is homogeneous on the zero section $G/H$ and hence agrees with a hyperholomorphic bundle there. Invoking \cite{Feix1},\cite{Feix2} once more we see that they are isomorphic everywhere.

\vskip 1cm
 Mathematical Institute, 24-29 St Giles, Oxford OX1 3LB, UK
 
 hitchin@maths.ox.ac.uk

 \end{document}